

\baselineskip=14pt
\parskip=10pt
\def\halmos{\hbox{\vrule height0.15cm width0.01cm\vbox{\hrule height
  0.01cm width0.2cm \vskip0.15cm \hrule height 0.01cm width0.2cm}\vrule
  height0.15cm width 0.01cm}}
\font\eightrm=cmr8 
\font\eighttt=cmtt8
\magnification=\magstephalf

\def\1{{\overline{1}}}
\def\2{{\overline{2}}}
\parindent=0pt
\overfullrule=0in

\def\frac#1#2{{#1 \over #2}}
\nopagenumbers
\bf
\centerline
{
A Maple One-Line Proof of George Andrews's  Formula that Says that the Number 
}
\centerline
{
of Triangles with Integer Sides Whose Perimeter is $n$ Equals $\{ \frac{n^2}{12} \}-\lfloor \frac{n}{4} \rfloor \lfloor \frac{n+2}{4} \rfloor$
}
\rm
\bigskip
\centerline{ {\it Shalosh B. EKHAD }\footnote{$^1$}
{\eightrm  \raggedright
{\eighttt ShaloshBEkhad@gmail.com} \quad .
\hfill \break
Exclusively published in: {\eighttt http://www.math.rutgers.edu/\~{}zeilberg/pj.html} and
{\eighttt http://arxiv.org/} . \hfill \break
Feb. 4 2012.
}
}

{\tt evalb(seq(coeff(taylor(q\^{}3/(1-q\^{}2)/(1-q\^{}3)/(1-q\^{}4),q=0,37),q,i),i=0..36) \hfill\break
=seq(round(n\^{}2/12)-trunc(n/4)*trunc((n+2)/4),n=0..36)); }

{\bf Comments by Doron Zeilberger}

{\bf 1}. The succinct formula of the title was discovered and first proved by George Andrews[A1],
in a less-than-one-page cute note, improving a four-page note by  Jordan et. al. [JWW].
Andrews's note, while short and sweet, is not self-contained, using results from his classic [A2].

{\bf 2}. Here is a clarification of the above Maple one-line proof. Arrange the lengths of the sides of a typical integer-side 
triangle in non-increasing order, and write them
as $[a+b+c+1,b+c+1,c+1]$ for $a,b,c \geq 0$. By [E] I.20,
$(b+c+1)+(c+1) >a+b+c+1$, so $c \geq a$, so $c=a+t$ for $t \geq 0$. So a generic integer-side triangle can be written
as $[a+b+(a+t)+1,b+(a+t)+1,a+t+1]=[2a+b+t+1,a+b+t+1,a+t+1]$, and hence the number of triangles with integer sides whose
perimeter equals $n$ is the coefficient of $q^n$ in 
$$
\sum_{a,b,t \geq 0} q^{(2a+b+t+1)+(a+b+t+1)+(a+t+1)}=\sum_{a,b,t \geq 0} q^{4a+3t+2b+3}
=q^3 \left ( \sum_{a \geq 0} q^{4a} \right ) \left ( \sum_{t \geq 0} q^{3t} \right ) \left ( \sum_{b \geq 0} q^{2b} \right )=
$$
$$
\frac{q^3}{(1-q^4)(1-q^3)(1-q^2)} \quad .
$$
The coefficient of $q^n$ is obviously a {\it quasi-polynomial} of degree $2$ and ``period'' $lcm (2,3,4)=12$, but so is 
$\{ \frac{n^2}{12} \}-\lfloor \frac{n}{4} \rfloor \lfloor \frac{n+2}{4} \rfloor$, hence it suffices to check
that the first $(2+1) \cdot 12+1=37$ values match. \halmos

{\bf 3}. This is yet another example where ``physical'' (as opposed to ``mathematical'') induction, i.e. checking
{\bf finitely} (and not that many!) special cases, constitutes a {\it perfectly rigorous proof}! No need for
the ``simple'' argument by ``mathematical'' induction alluded to by Andrews in his note. In this case we were
in the {\it quasi-polynomial} {\bf ansatz}. See [Z1][Z2], as well as the future classic [KP].

{\bf References}

[A1] G. E. Andrews, {\it A note on partitions and triangles with integer sides}, 
Amer. Math. Monthly {\bf 86} (1979),477-478.

[A2] G. E. Andrews, {\it ``The Theory of Partitions''}, Addison-Wesley, 1976. Reprinted by 
Cambridge University Press, 1984. First paperback edition, 1998.

[E] Euclid, {\it ``The Elements''}, Alexandria University Press, ca. 300 BC.

[JWW] J.H. Jordan, R. Walch, and R.J. Wisner, {\it Triangles with integer sides},
Amer. Math. Monthly {\bf 86} (1979), 686-689.

[KP] M. Kauers and P. Paule, {\it ``The Concrete Tetrahedron''}, Springer, 2011 \hfill\break
{\tt http://www.springer.com/mathematics/analysis/book/978-3-7091-0444-6} \quad .

[Z1] D. Zeilberger, {\it Enumerative and Algebraic Combinatorics},
in: ``Princeton Companion to Mathematics'' , (Timothy Gowers, ed.), Princeton University Press, pp. 550-561, 2008. \hfill\break
{\tt http://www.math.rutgers.edu/\~{}zeilberg/mamarim/mamarimPDF/enu.pdf} 

[Z2] D. Zeilberger, {\it An Enquiry Concerning Human (and Computer!) [Mathematical] Understanding},
in: C.S. Calude, ed., ``Randomness \& Complexity, from Leibniz to Chaitin'' World Scientific, Singapore, 2007. \hfill\break
{\tt http://www.math.rutgers.edu/\~{}zeilberg/mamarim/mamarimPDF/enquiry.pdf}
\end